\documentclass[11pt]{article}

\usepackage[letterpaper,margin=1in]{geometry}
\usepackage[T1]{fontenc}
\usepackage[utf8]{inputenc}

\usepackage{amsmath}
\usepackage{amssymb}
\usepackage{amsthm}
\usepackage{mathtools}

\usepackage[defaultsups]{newtxtext}
\usepackage{newtxmath}
\usepackage{bm}
\usepackage{bbm}

\usepackage{graphicx}
\usepackage[shortlabels]{enumitem}
\usepackage{microtype}
\usepackage{booktabs,tabularx,array,longtable}
\usepackage{authblk}   % must precede hyperref

\allowdisplaybreaks

\usepackage[numbers,sort&compress]{natbib}
\bibpunct[, ]{[}{]}{,}{n}{}{,}%
\usepackage{xcolor}
\definecolor{linkcolour}{rgb}{0.05,0.25,0.55}
\usepackage[%
  colorlinks=true,
  linkcolor=linkcolour,
  citecolor=linkcolour,
  urlcolor=linkcolour,
    pdftitle={Optimal Last-Iterate Bounds for Nonexpansive Stochastic Approximation},
  pdfkeywords={},
  pdfauthor={Ramanan Tamizholi, Siddharth Chandak, Rahul Singh}%
]{hyperref}

\theoremstyle{plain}
\newtheorem{theorem}{Theorem}
\newtheorem{lemma}{Lemma}
\newtheorem{proposition}{Proposition}

\newtheorem{assumption}{Assumption}

\theoremstyle{definition}
\newtheorem{example}{Example}

\newtheorem{remark}{Remark}

\newcommand{\Xstar}{\mathcal{X}^\star}

\newcommand{\EE}{\mathbb{E}}
\newcommand{\cfrak}{\mathfrak{c}}

\newcolumntype{Y}{>{\raggedright\arraybackslash}X}

\title{\bfseries Optimal Last-Iterate Bounds for Nonexpansive Stochastic Approximation}

\author[1]{Ramanan Tamizholi}
\author[2]{Siddharth Chandak}
\author[3]{Rahul Singh}

\affil[1]{Department of Computer Science and Automation, Indian Institute of Science, Bengaluru, India. \texttt{ramanant@iisc.ac.in}}
\affil[2]{Department of Electrical Engineering, Stanford University, Stanford, California, USA. \texttt{chandaks@stanford.edu}}
\affil[3]{EPITA, Laboratoire Recherche de l'EPITA, Paris, France. \texttt{rahulsingh0188@gmail.com}}

\date{}

\begin{document}

\maketitle

\begin{abstract}
We establish last-iterate mean-square residual bounds for nonexpansive stochastic approximation (SA), or stochastic Krasnoselskii--Mann (KM) iterations. In particular, we consider fixed-point mappings that are nonexpansive with respect to the Euclidean norm, and focus on polynomially decaying stepsizes. We study three conditions on the conditional second moment of the noise, and obtain bounds that are optimal up to logarithmic factors in each case. When it scales linearly with the squared residual, the best rate is $O(1/k^{1-\epsilon})$ for any $\epsilon>0$. In this case, for all considered stepsizes, our bounds match those of the noiseless KM iteration. When it is bounded by a constant, the best rate is $O(\log^2 k/\sqrt{k})$. When it scales affinely with the squared norm of the iterate, the best rate is $O(1/k^{1/2-\epsilon})$, and the mean-square residual can diverge if the stepsizes are not square-summable. Thus, unlike for contractive SA, square-summable stepsizes are necessary for nonexpansive SA. We demonstrate how various algorithms with cocoercive maps, such as stochastic gradient descent, linear SA, and learning in cocoercive games, fit into our framework, and how the different noise conditions arise in each of them.
\end{abstract}

\smallskip

\section{Introduction}\label{sec:introduction}
Stochastic approximation (SA) is a class of iterative algorithms to find fixed points of an operator given its noisy realizations~\cite{RobbinsMonro1951}. These algorithms have found applications in various areas, including optimization, reinforcement learning~\cite{SuttonBarto2018}, communications and stochastic control~\cite{Borkar2022}. Consider the SA iteration
\begin{align*}
x_{k+1}=x_k+\beta_k\left(f(x_k)-x_k+M_{k+1}\right).
\end{align*}
Here, $x_k\in\mathbb{R}^d$ is the iterate, $\beta_k$ is the stepsize, $\{M_{k+1}\}$ is a martingale difference noise sequence, and $f:\mathbb{R}^d\to\mathbb{R}^d$ is the mapping for which we wish to find a fixed point $x^\star$, i.e., a point satisfying $f(x^\star)=x^\star$.

Most finite-time analyses of SA assume that $f$ is a contraction~\cite{SrikantYing2019,ChenEtAl2020,ChandakBorkarDodhia2022,Chandak2026UnifiedContractive}. However, this assumption does not hold in many problems of interest, e.g., stochastic gradient descent (SGD) for convex functions that are not strongly convex, where the gradient step is only nonexpansive. In this work, we study SA where $f$ is nonexpansive with respect to the Euclidean norm. These iterations are also known as stochastic Krasnoselskii--Mann (KM) iterations~\cite{Krasnoselskii1955,Mann1953}. As nonexpansive mappings need not have a unique fixed point, we bound the last-iterate mean-square residual $\EE\left[\|f(x_k)-x_k\|^2\right]$ rather than the distance to a fixed point.

While there has been recent interest in stochastic fixed-point computation for nonexpansive mappings, existing works differ from ours in two ways. First, some works consider mappings that are nonexpansive with respect to general norms~\cite{BravoCominetti2024,BlaserZhang2026}. Due to their more general scope, their bounds are weaker than what can be obtained in the Euclidean setting. Second, other works modify the iteration, either through Halpern-type anchoring~\cite{CaiEtAl2022,BravoContreras2026,Diakonikolas2026,Shehu2026} or through modifications of the KM iteration such as minibatches or inertia~\cite{Iiduka2026,YangEtAl2026}. In contrast, we focus on the standard stochastic KM iteration with polynomially decaying stepsizes under the Euclidean norm.

We study three conditions on the conditional second moment of the noise. Under each condition, we obtain the first last-iterate mean-square residual bounds that are optimal up to logarithmic factors.
\begin{itemize}
\item[\textup{($\mathrm R$)}] \textbf{R}esidual-dependent: bounded by a multiple of the squared residual $\|f(x_k)-x_k\|^2$.
\begin{itemize}
\item Best rate: $O(1/k^{1-\epsilon})$, where $\epsilon>0$ can be arbitrarily small.
\end{itemize}
\item[\textup{($\mathrm C$)}] Bounded by a \textbf{C}onstant.
\begin{itemize}
\item Best rate: $\widetilde O(1/\sqrt{k})$, where $\widetilde O$ hides logarithmic factors.
\end{itemize}
\item[\textup{($\mathrm I$)}] \textbf{I}terate-dependent: bounded by an affine function of the squared norm of the iterate $\|x_k\|^2$.
\begin{itemize}
\item Best rate: $O(1/k^{1/2-\epsilon})$, where $\epsilon>0$ can be arbitrarily small.
\end{itemize}
\end{itemize}
While most finite-time guarantees for stochastic KM iterations are obtained under condition $\mathrm C$~\cite{BravoCominetti2024,BravoContreras2026}, condition $\mathrm I$ is the standard noise model in the SA literature~\cite{Borkar2022}. We additionally include condition $\mathrm R$ for two reasons. First, it arises in applications, e.g., as the strong growth condition for SGD~\cite{SchmidtLeRoux2013,VaswaniBachSchmidt2019}. Second, it shows that noise vanishing at the fixed points does not slow down convergence: under condition $\mathrm R$, our bounds match those for the noiseless KM iteration for all polynomially decaying stepsizes, which is not the case under condition $\mathrm C$ or $\mathrm I$. Moreover, under condition $\mathrm I$, we show that square-summable stepsizes are necessary, as the mean-square residual can otherwise diverge. This requirement is absent for contractive SA, where convergence is achieved even for non-square-summable stepsizes~\cite{ChenEtAl2020,Chandak2026UnifiedContractive}, and thus just moving from contractive to nonexpansive mappings brings it back. We explain the relevance of these conditions through three applications: SGD, linear SA, and learning in games.

Finally, we briefly outline our proof technique. Bounding the ergodic residual, i.e., a weighted average of past residuals, is well understood from the analysis of inexact KM iterations~\cite{DavisYin2016,LiangFadiliPeyre2016}. To obtain the last-iterate bound, we compare the stochastic iterates with a noiseless auxiliary iteration restarted from an earlier stochastic iterate, and bound how far the two trajectories can drift apart.

\subsection{Related Work}\label{subsec:related-work}
In recent years, stochastic KM iterations have been studied for mappings that are nonexpansive with respect to general norms. In this setting, almost sure convergence and finite-time last-iterate expected residual bounds were established in~\cite{BravoCominetti2024}, and extended to Markovian noise in~\cite{BlaserZhang2026}. In contrast, we focus on the Euclidean norm, which allows us to obtain sharper bounds on the mean-square residual. Moreover, while the analyses in~\cite{BravoCominetti2024,BlaserZhang2026} also rely on an auxiliary iteration, it is defined through an averaged noise sequence rather than a noiseless iteration.

Much of the remaining literature on nonexpansive mappings either modifies the KM iteration or studies variants of Halpern iterations. The former includes using increasing batch sizes~\cite{Iiduka2026} and introducing inertia~\cite{YangEtAl2026}, and the latter includes Halpern iterations combined with minibatches~\cite{BravoContreras2026,IidukaHalpern2026}, variance reduction~\cite{CaiEtAl2022,Diakonikolas2026,DiakonikolasGuzmanMartinezRubio2026,KayaEtAl2026}, or tuned anchoring schedules~\cite{Shehu2026}.

Bounds on the convergence rates of two-time-scale SA with a nonexpansive mapping in the slower time-scale have also been obtained in~\cite{Chandak2026Nonexpansive,SarkarAggarwal2026}. While the slower time-scale in this setting is also a stochastic inexact KM iteration, its rates are further limited by the time-scale separation.

Closest to our work, \cite{CortildCartisPeypouquet2026} studies stochastic KM iterations in Hilbert spaces. It assumes that each sampled function is nonexpansive and that the noise has finite variance at a fixed point, and establishes almost sure convergence and finite-time last-iterate mean-square residual bounds. In contrast, we require only the underlying function $f$ to be nonexpansive and impose moment conditions directly on the noise, without requiring it to arise from sampled nonexpansive functions. Moreover, our bounds improve on those of~\cite{CortildCartisPeypouquet2026} for some stepsizes and match them for others.

\subsection{Outline and Notation}\label{subsec:outline-notation}
The paper is organized as follows. Section~\ref{sec:main} sets up the problem and presents our upper bounds under each noise condition, followed by remarks on the choice of stepsizes and the optimality of the bounds. Section~\ref{sec:outline} gives a proof sketch of the upper bounds. Section~\ref{sec:lower} presents the lower bounds, which establish optimality, and the divergence result under condition $\mathrm I$. Section~\ref{sec:applications} discusses three applications and how our noise conditions arise in them. Finally, Section~\ref{sec:conclusion} concludes the paper.

Throughout, $\|\cdot\|$ and $\langle\cdot,\cdot\rangle$ denote the Euclidean norm and inner product on $\mathbb{R}^d$. For sequences, $O(\cdot)$ and $\Omega(\cdot)$ denote the standard order notation, and $\widetilde{O}(\cdot)$ additionally hides logarithmic factors. 
\section{Problem Formulation and Main Results}\label{sec:main}
Consider the following iteration.
\begin{align}
x_{k+1}=x_k+\beta_k\left(f(x_k)-x_k+M_{k+1}\right).\label{eq:iteration}
\end{align}
Here, $x_k\in\mathbb{R}^d$ is the iterate and $\beta_k$ is the stepsize at time $k\ge0$. The function $f:\mathbb{R}^d\to\mathbb{R}^d$ denotes the mapping whose fixed point we wish to find, and $M_{k+1}$ denotes the martingale difference noise sequence. 

The first assumption states that the mapping $f$ is nonexpansive and has at least one fixed point.
\begin{assumption}\label{assm:nonexp}
The function $f$ is nonexpansive, i.e., for all $x_1,x_2\in\mathbb{R}^d$,
\begin{align*}
\|f(x_1)-f(x_2)\|\le\|x_1-x_2\|.
\end{align*}
Moreover, the set of fixed points $\Xstar=\{x\in\mathbb{R}^d:f(x)=x\}$ is nonempty.
\end{assumption}
Unlike a contraction, a nonexpansive mapping can have multiple fixed points. Hence, instead of the distance $\|x_k-x^\star\|$ to a particular fixed point $x^\star$, we measure convergence through the residual $\|f(x_k)-x_k\|$, which is the standard measure of convergence for nonexpansive SA and KM iterations. Our objective is to obtain finite-time bounds on the last-iterate mean-square residual $\EE\left[\|f(x_k)-x_k\|^2\right]$.

We next state the assumption on the noise sequence $\{M_{k+1}\}$. We consider three different conditions on its conditional second moment, and obtain a separate bound under each of them.
\begin{assumption}\label{assm:noise}
Define the family of $\sigma$-fields $\mathcal{F}_k=\sigma(x_0,M_i,i\le k)$. Then $\{M_{k+1}\}$ is a martingale difference sequence with respect to $\mathcal{F}_k$, i.e., $\EE\left[M_{k+1}\mid\mathcal{F}_k\right]=0$. Moreover, for some constant $\sigma>0$, one of the following conditions holds for all $k\ge0$.
\begin{enumerate}
\item[\textup{($\mathrm R$)}] Conditional second moment is \textbf{R}esidual-dependent:
\begin{align*}
\EE\left[\|M_{k+1}\|^2\mid\mathcal{F}_k\right]\le\sigma^2\|f(x_k)-x_k\|^2.
\end{align*}
\item[\textup{($\mathrm C$)}] Conditional second moment is bounded by a \textbf{C}onstant:
\begin{align*}
\EE\left[\|M_{k+1}\|^2\mid\mathcal{F}_k\right]\le\sigma^2.
\end{align*}
\item[\textup{($\mathrm I$)}] Conditional second moment is affinely bounded by the \textbf{I}terate:
\begin{align*}
\EE\left[\|M_{k+1}\|^2\mid\mathcal{F}_k\right]\le\sigma^2\left(1+\|x_k\|^2\right).
\end{align*}
\end{enumerate}
\end{assumption}
Under condition $\mathrm R$, the noise scales with the residual, and hence vanishes at the fixed points of $f$. For example, for SGD iterations, noise whose conditional second moment scales with the squared norm of the gradient satisfies condition $\mathrm R$~\cite{BottouCurtisNocedal2018,SchmidtLeRoux2013,VaswaniBachSchmidt2019}. Condition $\mathrm I$ is the most general of the three, and is the standard noise model in the SA literature~\cite{Borkar2022,KushnerYin2003}. It clearly subsumes condition $\mathrm C$, and also subsumes condition $\mathrm R$, albeit with $\sigma^2$ replaced by $8\sigma^2\max\{1,\|x^\star\|^2\}$ for any $x^\star\in\Xstar$. However, the more restrictive conditions lead to better rates, and we therefore treat the three conditions separately.

Finally, we state the assumption on the stepsize sequence.
\begin{assumption}\label{assm:stepsize}
The stepsize sequence is of the form
\begin{align*}
\beta_k=\frac{\beta}{(k+K_0)^b},
\end{align*}
where $b\in(0,1)$, $0<\beta\le1/2$, and $K_0\ge1$. Moreover, under condition $\mathrm R$ in Assumption~\ref{assm:noise}, $(1+\sigma^2)\beta_0<1$.
\end{assumption}
The conditions $\beta\le1/2$ and $K_0\ge1$ ensure that $\beta_k\le1/2$ for all $k\ge0$. The additional condition under noise condition $\mathrm R$ ensures that the stepsize is sufficiently small relative to the noise magnitude, starting from time $k=0$. We do not consider stepsizes of the form $\beta_k=\Theta(1/k)$, as they lead to bounds of the form $O(1/\log k)$, which are weaker than the bounds obtained for $b<1$.

We now state our main results, one for each noise condition in Assumption~\ref{assm:noise}. An outline of their proofs is presented in Section~\ref{sec:outline}, and the complete proofs are given in Appendix~\ref{app:restart-proofs}. The lower bounds showing optimality of these results are presented in Section~\ref{sec:lower}. Each theorem is followed by a short discussion of the resulting rate and the choice of stepsize.
\begin{theorem}\label{thm:main_R}
Suppose Assumptions~\ref{assm:nonexp}--\ref{assm:stepsize} hold, with noise condition $\mathrm R$ in Assumption~\ref{assm:noise}.
Then there exists a constant $\cfrak_1>0$ such that, for all $k\geq0$,
\begin{align*}
\EE\left[\|f(x_k)-x_k\|^2\right] &\leq \dfrac{\cfrak_1}{(k+1)^{1-b}}.
\end{align*}
Moreover, this bound is optimal (see Remark~\ref{rem:optimal}).
\end{theorem}
\paragraph{Best rate.} A slowly decreasing stepsize, i.e., a smaller $b$, gives a better rate. In fact, we achieve the same rate as that achieved for noiseless KM iterations~\cite{CominettiSotoVaisman2014,DavisYin2016}. The best rate is $O(1/k^{1-\epsilon})$, where $\epsilon>0$ can be arbitrarily small, and is achieved by choosing $b=\epsilon$.
\begin{theorem}\label{thm:main_C}
Suppose Assumptions~\ref{assm:nonexp}--\ref{assm:stepsize} hold, with noise condition $\mathrm C$ in Assumption~\ref{assm:noise}.
Then there exist constants $\cfrak_2,\cfrak_3>0$ such that, for all $k\geq0$,
\begin{align*}
\EE\left[\|f(x_k)-x_k\|^2\right] &\leq \begin{cases}
\dfrac{\cfrak_2\log^2(k+2)}{(k+1)^b},&0<b\le\tfrac12,\\[8pt]
\dfrac{\cfrak_3}{(k+1)^{1-b}},&\tfrac12<b<1.
\end{cases}
\end{align*}
Moreover, this bound is optimal up to a logarithmic factor (see Remark~\ref{rem:optimal}).
\end{theorem}
\paragraph{Non-monotone dependence on the stepsize.} Unlike in Theorem~\ref{thm:main_R}, the rate is no longer monotone in $b$. For $b>1/2$, a smaller $b$ gives a better rate, as before. For $b<1/2$, the noise dominates, and the bound decays only as $(k+1)^{-b}$, up to logarithmic factors. Hence, the best rate of $\widetilde{O}(1/\sqrt{k})$ is achieved by choosing $b=1/2$.
\begin{theorem}\label{thm:main_I}
Suppose Assumptions~\ref{assm:nonexp}--\ref{assm:stepsize} hold, with noise condition $\mathrm I$ in Assumption~\ref{assm:noise}, and let $\tfrac12<b<1$.
Then there exists a constant $\cfrak_4>0$ such that, for all $k\geq0$,
\begin{align*}
\EE\left[\|f(x_k)-x_k\|^2\right] &\leq \dfrac{\cfrak_4}{(k+1)^{1-b}}.
\end{align*}
Moreover, this bound is optimal (see Remark~\ref{rem:optimal}). However, this result cannot be extended to $b\le1/2$: for every $0<b\le1/2$, there exists an instance satisfying Assumptions~\ref{assm:nonexp}--\ref{assm:stepsize} with noise condition $\mathrm I$ for which $\EE\left[\|f(x_k)-x_k\|^2\right]\to\infty$ as $k\to\infty$ (see Proposition~\ref{prop:copy-lower-I}).
\end{theorem}
\paragraph{Stepsize choice and necessity of square-summable stepsizes.} For $b>1/2$, the rate matches that of Theorems~\ref{thm:main_R} and~\ref{thm:main_C}. Since this result cannot be extended to $b\le1/2$, the best rate is $O(1/k^{1/2-\epsilon})$, where $\epsilon>0$ can be arbitrarily small, and is achieved by choosing $b=1/2+\epsilon$. The condition $b>1/2$ is equivalent to the stepsizes being square-summable. While this is the standard requirement for the convergence analysis of SA without additional assumptions on the operator~\cite{RobbinsMonro1951,Borkar2022,KushnerYin2003}, it has been shown that for contractive fixed-point iterations, convergence is achieved for all $0<b\le1$~\cite{ChenEtAl2020,Chandak2026UnifiedContractive}. Theorem~\ref{thm:main_I} shows that just by moving from contractive to nonexpansive mappings, the requirement of square-summable stepsizes comes back. This is not an artifact of our analysis: Proposition~\ref{prop:copy-lower-I} gives a single instance for which the mean-square residual grows as $\exp\left(\Omega(k^{1-2b})\right)$ for $0<b<1/2$, and polynomially for $b=1/2$. In contrast, Theorems~\ref{thm:main_R} and~\ref{thm:main_C} hold for all $0<b<1$, and do not require square-summable stepsizes. In fact, under noise condition $\mathrm C$, the best rate is achieved by $b=1/2$.

The following remarks discuss the optimality of the above bounds, and the implied bounds on the ergodic residual.
\begin{remark}[Optimality of the bounds]\label{rem:optimal}
By optimality in Theorems~\ref{thm:main_R}--\ref{thm:main_I}, we mean the following. For every $k\ge1$, there exist a mapping $f$ satisfying Assumption~\ref{assm:nonexp}, an initial iterate $x_0$, and a noise sequence satisfying the corresponding noise condition, for which $\EE\left[\|f(x_k)-x_k\|^2\right]$ is lower bounded by the upper bound at time $k$, up to a constant independent of $k$ (and up to a logarithmic factor under noise condition $\mathrm C$ with $b\le1/2$). Here, the instance is allowed to depend on $k$. The corresponding instances are given in Propositions~\ref{prop:copy-lower-R} and~\ref{prop:copy-lower-C}.
\end{remark}
\begin{remark}[Ergodic bounds]\label{rem:ergodic}
Another object of interest when studying nonexpansive SA is the ergodic residual $\frac{1}{k+1}\sum_{i=0}^{k}\EE\left[\|f(x_i)-x_i\|^2\right]$, or the minimum residual $\min_{0\le i\le k}\EE\left[\|f(x_i)-x_i\|^2\right]$~\cite{DavisYin2016}. We do not highlight bounds on these quantities separately here, because they follow directly from our last-iterate bounds. Both quantities are bounded by the average of the last-iterate bounds in Theorems~\ref{thm:main_R}--\ref{thm:main_I}, and hence satisfy the same upper bounds, up to constants.
\end{remark}
\section{Proof Outline for the Upper Bounds}\label{sec:outline}
In this section, we present an outline of the proofs of the upper bounds in Theorems~\ref{thm:main_R}--\ref{thm:main_I} through a series of lemmas. The complete proofs are given in Appendix~\ref{app:restart-proofs}. Throughout this section, we fix a fixed point $x^\star\in\Xstar$, and for brevity, we write $r_i:=\EE\left[\|f(x_i)-x_i\|^2\right]$ for the mean-square residual at time $i$, and
\begin{align*}
\bar r_{\ell,k}:=\frac{\sum_{i=\ell}^k\beta_ir_i}{\sum_{i=\ell}^k\beta_i},\qquad 0\le\ell\le k,
\end{align*}
for the weighted average of the residuals between times $\ell$ and $k$. Note that $\bar r_{k,k}=r_k$ is the last-iterate residual.

The proof consists of two steps. In the first step, we bound the ergodic residual $\bar r_{0,k}$. This follows from a standard Lyapunov drift argument for the squared distance $\|x_k-x^\star\|^2$, although some care is needed in handling the noise under each noise condition. In the second step, we convert this ergodic bound into a last-iterate bound on $r_k$. To do so, we compare the stochastic iterates with a noiseless trajectory restarted from an intermediate iterate. The residual along a noiseless KM trajectory is nonincreasing, and we show that the residuals of the two trajectories remain close in a suitable averaged sense.

The majority of the proof is common to the three noise conditions. The noise condition only enters through the bound on the accumulated noise $\sum_i\beta_i^2\EE\left[\|M_{i+1}\|^2\right]$, which appears in both steps.

\subsection{Ergodic Bound}\label{sec:outline-ergodic}
The first lemma is a one-step drift inequality for the squared distance of the iterates from $x^\star$. It follows from the identity $\|(1-\beta)a+\beta b\|^2=(1-\beta)\|a\|^2+\beta\|b\|^2-\beta(1-\beta)\|a-b\|^2$, the nonexpansiveness of $f$, and the martingale difference property of the noise.
\begin{lemma}\label{lem:lyapunov}
For every $k\ge0$,
\begin{align}
\EE\left[\|x_{k+1}-x^\star\|^2\right] &\le\EE\left[\|x_k-x^\star\|^2\right]-\beta_k(1-\beta_k)\EE\left[\|f(x_k)-x_k\|^2\right]+\beta_k^2\EE\left[\|M_{k+1}\|^2\right].\label{eq:lyapunov}
\end{align}
\end{lemma}
When summed over $k$, the squared distance terms telescope, and dividing by $\sum_{i=0}^k\beta_i$ gives a bound on the ergodic residual $\bar r_{0,k}$ in terms of $\|x_0-x^\star\|^2$ and the accumulated noise $\sum_{i=0}^k\beta_i^2\EE\left[\|M_{i+1}\|^2\right]$. The three noise conditions differ only in how the accumulated noise is handled.
\begin{itemize}
\item Under noise condition $\mathrm R$, the noise term $\sigma^2\beta_i^2r_i$ is absorbed into the drift term $-\beta_i(1-\beta_i)r_i$, which is possible since $(1+\sigma^2)\beta_i\le(1+\sigma^2)\beta_0<1$ by Assumption~\ref{assm:stepsize}. The resulting ergodic bound is of the order of $1/\sum_{i=0}^k\beta_i$, the same as for the noiseless KM iteration.
\item Under noise condition $\mathrm C$, the noise contributes the additional term $\sigma^2\sum_{i=0}^k\beta_i^2$. For $b>1/2$, this sum is bounded, and the ergodic bound is of the order of $1/\sum_{i=0}^k\beta_i$, as under noise condition $\mathrm R$. For $b<1/2$, it grows as $(k+1)^{1-2b}$ and dominates, giving an ergodic bound of the order of $(k+1)^{-b}$, i.e., of the order of the stepsize $\beta_k$.
\item Under noise condition $\mathrm I$, we first show that for $b>1/2$, the iterates are bounded in expectation, i.e., $\EE\left[\|x_k-x^\star\|^2\right]$ is bounded uniformly in $k$. This is the only step where $b>1/2$ is needed, and it cannot be extended to $b\le1/2$. Indeed, since $\|f(x)-x\|\le2\|x-x^\star\|$, a uniform bound on $\EE\left[\|x_k-x^\star\|^2\right]$ would also give a uniform bound on the residual $r_k$, while Proposition~\ref{prop:copy-lower-I} shows that $r_k$ can diverge for $b\le1/2$. Then, the noise second moment is bounded by a constant, and the argument reduces to that under noise condition $\mathrm C$.
\end{itemize}
\begin{lemma}\label{lem:residual-sum}
For every $k\ge0$, the following hold.
\begin{enumerate}
\item[(a)] Under noise condition $\mathrm R$,
\begin{align}
\bar r_{0,k} &\le\frac{\|x_0-x^\star\|^2}{[1-(1+\sigma^2)\beta_0]\sum_{i=0}^k\beta_i}.\label{eq:pure-sum}
\end{align}
\item[(b)] Under noise condition $\mathrm C$,
\begin{align}
\bar r_{0,k} &\le\frac{2\|x_0-x^\star\|^2+2\sigma^2\sum_{i=0}^k\beta_i^2}{\sum_{i=0}^k\beta_i}.\label{eq:C-sum}
\end{align}
\item[(c)] Let $\Gamma_1:=\sum_{i=0}^\infty\beta_i^2$. Under noise condition $\mathrm I$ with $1/2<b<1$, there exists a constant $\Gamma_2>0$, independent of $k$, such that $\EE\left[\|x_k-x^\star\|^2\right]\le\Gamma_2$, and
\begin{align}
\bar r_{0,k} &\le\frac{2\left(\|x_0-x^\star\|^2+\sigma^2\Gamma_1(1+2\|x^\star\|^2+2\Gamma_2)\right)}{\sum_{i=0}^k\beta_i}.\label{eq:I-finite}
\end{align}
\end{enumerate}
\end{lemma}

\subsection{Last-Iterate Bound}\label{sec:outline-last}
The key step in going from the ergodic residual to the last iterate is an auxiliary noiseless iteration, restarted at an arbitrary iterate of the stochastic trajectory. As discussed above, the residual along a noiseless trajectory is nonincreasing. We use this auxiliary iteration to show that the residuals are approximately nonincreasing on average, in the sense that $\bar r_{\ell+1,k}$ exceeds $\bar r_{\ell,k}$ by at most a noise-dependent term. For a fixed $\ell\ge0$, define the auxiliary iteration
\begin{align}\label{eq:copy}
y_\ell=x_\ell,\qquad y_{k+1}=y_k+\beta_k\left(f(y_k)-y_k\right),\quad k\ge \ell.
\end{align}
Since the two trajectories agree at time $\ell$, any discrepancy between them after time $\ell$ is due to the noise injected from time $\ell$ onward. The following lemma quantifies this discrepancy. Its first inequality is the analogue of Lemma~\ref{lem:lyapunov} for the difference $x_k-y_k$. In fact, as the constant sequence $y_k=x^\star$ is itself a noiseless trajectory, Lemma~\ref{lem:lyapunov} is the special case of this inequality with $y_k=x^\star$. Summing the first inequality, as before, gives the second inequality: the weighted sum of the squared differences between the residuals of the two trajectories is bounded by the noise injected after time $\ell$. The final statement shows that the residual along the noiseless trajectory is nonincreasing, and is a consequence of the nonexpansiveness of $f$.
\begin{lemma}\label{lem:copy-diss}
For every $k\ge \ell$,
\begin{equation}
\EE\left[\|x_{k+1}-y_{k+1}\|^2\right] \le\EE\left[\|x_k-y_k\|^2\right]
-\beta_k(1-\beta_k)\EE\left[\|f(x_k)-x_k-f(y_k)+y_k\|^2\right]+\beta_k^2\EE\left[\|M_{k+1}\|^2\right],\label{eq:copy-diss}\end{equation}
and 
\begin{equation}\sum_{i=\ell}^k\beta_i\EE\left[\|f(x_i)-x_i-f(y_i)+y_i\|^2\right] \le2\sum_{i=\ell}^k\beta_i^2\EE\left[\|M_{i+1}\|^2\right].\label{eq:copy-sum}
\end{equation}
Moreover,
\begin{align*}
\|f(y_k)-y_k\|\le\|f(y_\ell)-y_\ell\|=\|f(x_\ell)-x_\ell\|.
\end{align*}
\end{lemma}
To move from the weighted average $\bar r_{m,k}$ to the last-iterate residual $\bar r_{k,k}=r_k$, we need to bound how $\bar r_{\ell,k}$ grows as $\ell$ increases. Removing the first term from the weighted average gives the exact identity
\begin{align*}
\bar r_{\ell+1,k}=\bar r_{\ell,k}+\frac{\beta_\ell}{\sum_{i=\ell+1}^k\beta_i}\left(\bar r_{\ell,k}-r_\ell\right).
\end{align*}
In the absence of noise, the monotonicity of the residual gives $\bar r_{\ell,k}\le r_\ell$, and hence the second term is nonpositive. In the presence of noise, this term can be bounded using Lemma~\ref{lem:copy-diss}, which gives the following lemma.
\begin{lemma}\label{lem:copy-future}
For any $0\le \ell<k$,
\begin{align}
\bar r_{\ell+1,k}\le\bar r_{\ell,k}+2\sqrt2\cdot\frac{\beta_\ell}{\sum_{i=\ell+1}^k\beta_i}\left(\bar r_{\ell,k}\right)^{1/2}\left(\frac{\sum_{i=\ell}^k\beta_i^2\EE\left[\|M_{i+1}\|^2\right]}{\sum_{i=\ell}^k\beta_i}\right)^{1/2}.\label{eq:future}
\end{align}
\end{lemma}
Fix a starting time $0\le m\le k$, which we choose later. Applying Lemma~\ref{lem:copy-future} successively for $\ell=m,m+1,\ldots,k-1$ and accumulating the error terms gives a bound on $r_k$ in terms of $\bar r_{m,k}$. Under noise condition $\mathrm R$, the error term in Lemma~\ref{lem:copy-future} is at most a multiple of $\bar r_{\ell,k}$ itself, and hence the resulting bound takes the form of a product, as in part (a) of the following lemma.
\begin{lemma}\label{lem:copy-propagate}
For any $0\le m\le k$,
\begin{enumerate}
\item[(a)] Under noise condition $\mathrm R$,
\begin{align}
r_k &\le\bar r_{m,k}\prod_{i=m}^{k-1}\left[1+2\sqrt2\sigma\frac{\beta_i\sqrt{\beta_i}}{\sum_{j=i+1}^k\beta_j}\right].\label{eq:prop-pure}
\end{align}
\item[(b)] Under noise condition $\mathrm C$ with $0<b<1$, or $\mathrm I$ with $1/2<b<1$, there exists a constant $\Gamma_3>0$, independent of $m$ and $k$, such that
\begin{align}
r_k\le\Gamma_3\left(\bar r_{m,k}+\left(\sum_{i=m}^{k-1}\frac{\beta_i\sqrt{\beta_i}}{\sum_{j=i+1}^k\beta_j}\right)^2\right).\label{eq:prop-CI}
\end{align}
\end{enumerate}
\end{lemma}
The choice of the starting time $m$ involves a trade-off. On the one hand, Lemma~\ref{lem:copy-propagate} bounds $r_k$ in terms of $\bar r_{m,k}$, while Lemma~\ref{lem:residual-sum} only bounds $\bar r_{0,k}$. Hence, we need $\bar r_{m,k}$ to be of the same order as $\bar r_{0,k}$. Since $r_i\ge0$, we have $\bar r_{m,k}\le\bar r_{0,k}\sum_{i=0}^k\beta_i/\sum_{i=m}^k\beta_i$, and hence this holds as long as the ratio $\sum_{i=0}^k\beta_i/\sum_{i=m}^k\beta_i$ is bounded. This rules out $m$ close to $k$. On the other hand, the error terms in Lemma~\ref{lem:copy-propagate} involve the sum of the ratios $\beta_i^{3/2}/\sum_{j=i+1}^k\beta_j$ over $m\le i<k$, which only increases as $m$ decreases, since starting earlier adds terms with larger stepsizes. We take $m=\lceil k/2\rceil$, which balances the two, as shown in the following lemma.
\begin{lemma}\label{lem:terminal-sum}
Let $m=\lceil k/2\rceil$. There exist constants $\Gamma_4,\Gamma_5>0$ such that, for every $k\ge0$,
\begin{align}
\bar r_{m,k} \le\Gamma_4\bar r_{0,k},\quad\text{and}\quad
\sum_{i=m}^{k-1}\frac{\beta_i\sqrt{\beta_i}}{\sum_{j=i+1}^k\beta_j} \le\frac{\Gamma_5\log(k+2)}{(k+1)^{b/2}}.\label{eq:terminal}
\end{align}
\end{lemma}
Taking $m=\lceil k/2\rceil$ and combining Lemmas~\ref{lem:residual-sum}, \ref{lem:copy-propagate}, and~\ref{lem:terminal-sum} gives the upper bounds in Theorems~\ref{thm:main_R}--\ref{thm:main_I}. In particular, under noise condition $\mathrm R$, the product in Lemma~\ref{lem:copy-propagate}(a) is bounded by a constant, and hence the last-iterate residual is of the same order as the ergodic residual.

\section{Optimality of the Upper Bounds}\label{sec:lower}
In this section, we present the instances showing that the upper bounds in Theorems~\ref{thm:main_R}--\ref{thm:main_I} are optimal, and that Theorem~\ref{thm:main_I} cannot be extended to $b\le1/2$. For each result, we first describe the instance, and then discuss what it shows. Throughout this section, the stepsizes satisfy Assumption~\ref{assm:stepsize}, and the complete proofs are given in Appendix~\ref{app:lower-proofs}.

The first two instances are built from rotations. For $\theta\in(0,\pi]$, let
\begin{align*}
R_\theta=\begin{pmatrix}\cos\theta&-\sin\theta\\\sin\theta&\cos\theta\end{pmatrix}
\end{align*}
denote the rotation of $\mathbb{R}^2$ by the angle $\theta$. Since $R_\theta$ is orthogonal, the map $x\mapsto R_\theta x$ is nonexpansive, and since $\det(R_\theta-I)=2(1-\cos\theta)>0$, its unique fixed point is the origin. Moreover, for every $x\in\mathbb{R}^2$ and $\beta\in[0,1]$,
\begin{align}
\|R_\theta x-x\|^2=2(1-\cos\theta)\|x\|^2,\qquad\|(1-\beta)x+\beta R_\theta x\|^2=\left[1-2(1-\cos\theta)\beta(1-\beta)\right]\|x\|^2.\label{eq:copy-single-identity}
\end{align}
That is, the squared residual is a fixed fraction $2(1-\cos\theta)$ of the squared norm, and each noiseless KM step shrinks the squared norm by the factor $1-2(1-\cos\theta)\beta(1-\beta)$. A small angle therefore makes the residual small relative to the squared norm, but also makes the iterates contract slowly. Our instances exploit this trade-off by choosing the angle depending on the time horizon $k$.

\subsection{Lower Bound without Noise}\label{sec:lower-R}
\begin{example}[Slow rotation]\label{ex:rotation}
Fix $k\ge1$. Let $d=2$, $x_0=(1,0)$, $M_{i+1}=0$ for all $i\ge0$, and $f(x)=R_\theta x$, where $\theta$ is chosen such that
\begin{align*}
2(1-\cos\theta)=\frac{1}{1+2\sum_{j=0}^{k-1}\beta_j}.
\end{align*}
\end{example}
To see why this choice of the angle is natural, write $\gamma=2(1-\cos\theta)$ and $\tau=\sum_{j=0}^{k-1}\beta_j$. By \eqref{eq:copy-single-identity}, the squared residual at time $k$ is $\gamma\prod_{j=0}^{k-1}[1-\gamma\beta_j(1-\beta_j)]\|x_0\|^2$, which is roughly $\gamma e^{-\gamma\tau}\|x_0\|^2$. As a function of $\gamma$, this is maximized at $\gamma=1/\tau$, where it equals $\|x_0\|^2/(e\tau)$. The choice $\gamma=1/(1+2\tau)$ in Example~\ref{ex:rotation} is of this order, and hence gives a residual of the order of $1/\tau$. Since there is no noise, this instance satisfies all three noise conditions in Assumption~\ref{assm:noise}. Example~\ref{ex:rotation} is essentially the instance used in~\cite{SarkarAggarwal2026} to lower bound the residual of noiseless KM iterations, and we include it here for completeness.
\begin{proposition}\label{prop:copy-lower-R}
There exists a constant $\cfrak_5>0$ such that, for every $k\ge1$, the instance in Example~\ref{ex:rotation} satisfies
\begin{align*}
\EE\left[\|f(x_k)-x_k\|^2\right]\ge\frac{\cfrak_5}{(k+1)^{1-b}}.
\end{align*}
\end{proposition}
Proposition~\ref{prop:copy-lower-R} shows that the upper bounds in Theorem~\ref{thm:main_R}, in Theorem~\ref{thm:main_C} for $b>1/2$, and in Theorem~\ref{thm:main_I} are optimal.

\subsection{Lower Bound for Noise Condition $\mathrm C$}\label{sec:lower-C}
\begin{example}[Rotations at multiple scales]\label{ex:multiscale}
Fix $k\ge1$. Let $d=2k$, and group the coordinates of $\mathbb{R}^{2k}$ into $k$ consecutive pairs, indexed by $p\in\{1,\ldots,k\}$. Let $x_0=0$, and let $f(x)=\operatorname{diag}(R_{\theta_1},\ldots,R_{\theta_k})x$ act as the rotation $R_{\theta_p}$ on the $p$-th pair, where $\theta_p$ is chosen such that
\begin{align*}
2(1-\cos\theta_p)=\frac{1}{1+2\sum_{j=p}^{k-1}\beta_j},\qquad1\le p\le k.
\end{align*}
Let $\{\varepsilon_i\}_{i\ge1}$ be independent random variables taking the values $1$ and $-1$ with probability $1/2$ each. Let $M_{i+1}$ have a single nonzero coordinate, equal to $\sigma\varepsilon_{i+1}$, in the first coordinate of the $(i+1)$-th pair for $0\le i<k$, and of the first pair for $i\ge k$.
\end{example}
Only the noise injected before time $k$ affects $x_k$. The noise injected from time $k$ onward is defined only so that $\{M_{i+1}\}$ is a noise sequence for all $i\ge0$. The map $f$ is nonexpansive with the origin as its unique fixed point, and the noise satisfies noise condition $\mathrm C$ with equality. The proof of the following proposition consists of two steps.
\begin{itemize}
\item \emph{Each pair is a copy of Example~\ref{ex:rotation}.} Before time $k$, the $p$-th pair receives noise only at time $p-1$, and is zero up to that time. Hence, from time $p$ onward, it evolves as a copy of Example~\ref{ex:rotation}, started with squared norm $\sigma^2\beta_{p-1}^2$, and with the angle $\theta_p$ tuned to the remaining horizon $\sum_{j=p}^{k-1}\beta_j$. As discussed after Example~\ref{ex:rotation}, this choice of the angle maximizes, up to a constant factor, the residual of the $p$-th pair at time $k$, which is then of the order of $\sigma^2\beta_{p-1}^2/(1+2\sum_{j=p}^{k-1}\beta_j)$. Since $f$ acts separately on each pair and the pairs are orthogonal, the squared residual at time $k$ is the sum of these contributions. That is, the instance attains the worst case for each noise injection simultaneously.
\item \emph{Evaluating the sum.} Since all terms in the sum are nonnegative, we can lower bound it by restricting it to the noise injected at times $k/2\le i<k$. The restricted sum is of the order of $\sigma^2\beta_k\log k$, i.e., of the order of $\sigma^2\log k/k^b$.
\end{itemize}
\begin{proposition}\label{prop:copy-lower-C}
Let $0<b\le1/2$. There exists a constant $\cfrak_6>0$ such that, for every $k\ge1$, the instance in Example~\ref{ex:multiscale} satisfies
\begin{align*}
\EE\left[\|f(x_k)-x_k\|^2\right]\ge\frac{\cfrak_6\sigma^2\log(k+1)}{(k+1)^b}.
\end{align*}
\end{proposition}
Proposition~\ref{prop:copy-lower-C} shows that the upper bound in Theorem~\ref{thm:main_C} for $b\le1/2$ is optimal up to a logarithmic factor.

\subsection{Divergence under Noise Condition $\mathrm I$}\label{sec:lower-I}
\begin{example}[Reflection with noise condition $\mathrm I$]\label{ex:reflection}
Let $d=2$, write $x=(u,v)\in\mathbb{R}^2$, and let $f(u,v)=(u,-v)$ and $x_0=(1,0)$. Writing $x_i=(u_i,v_i)$, let the two coordinates of $M_{i+1}$, conditionally on $\mathcal{F}_i$, independently take the values $\sigma u_i/\sqrt2$ and $-\sigma u_i/\sqrt2$ with probability $1/2$ each.
\end{example}
In this instance, $f$ is a reflection, and is hence nonexpansive. Its fixed points form the line $\{(u,0):u\in\mathbb{R}\}$, and $\|f(x)-x\|^2=4v^2$. The noise satisfies $\EE\left[\|M_{i+1}\|^2\mid\mathcal{F}_i\right]=\sigma^2u_i^2\le\sigma^2(1+\|x_i\|^2)$, and hence noise condition $\mathrm I$ holds. Along the line of fixed points, the iteration has no drift, and the $u$-coordinate evolves as $u_{i+1}=u_i+\beta_iM_{i+1}^{(1)}$. As the noise scales with $u_i$, this is a multiplicative random walk, with $\EE[u_{i+1}^2]=(1+\sigma^2\beta_i^2/2)\EE[u_i^2]$. Hence, $\EE[u_k^2]$ grows as $\exp(\Omega(\sum_{i=0}^{k-1}\beta_i^2))$, which diverges when $\sum_i\beta_i^2=\infty$. Since the noise in the $v$-coordinate also scales with $u_i$, this growth is transferred to the residual, up to a polynomial factor.
\begin{proposition}\label{prop:copy-lower-I}
The instance in Example~\ref{ex:reflection} satisfies the following.
\begin{enumerate}
\item[(a)] If $0<b<1/2$, there exist constants $\cfrak_7,\cfrak_8>0$ such that, for all $k\ge1$,
\begin{align*}
\EE\left[\|f(x_k)-x_k\|^2\right]\ge\frac{\cfrak_7}{(k+1)^{2b}}\exp\left(\cfrak_8(k+1)^{1-2b}\right).
\end{align*}
\item[(b)] If $b=1/2$ and $\sigma\beta>\sqrt2$, there exists a constant $\cfrak_9>0$ such that, for all $k\ge1$,
\begin{align*}
\EE\left[\|f(x_k)-x_k\|^2\right]\ge\cfrak_9(k+1)^{\sigma^2\beta^2/2-1}.
\end{align*}
\end{enumerate}
In particular, in both cases, $\EE\left[\|f(x_k)-x_k\|^2\right]\to\infty$ as $k\to\infty$.
\end{proposition}
Proposition~\ref{prop:copy-lower-I} shows that Theorem~\ref{thm:main_I} cannot be extended to $b\le1/2$. The condition $\sigma\beta>\sqrt2$ in part (b) is allowed, as noise condition $\mathrm I$ places no restriction on $\sigma$. Unlike Examples~\ref{ex:rotation} and~\ref{ex:multiscale}, this is a single instance that does not depend on $k$. As discussed above, the divergence is driven by the noise scaling with the iterate along the line of fixed points. Under noise condition $\mathrm C$, the $u$-coordinate would instead perform an additive random walk, whose second moment grows only as $\sum_{i=0}^{k-1}\beta_i^2$, and the residual would still be controlled by Theorem~\ref{thm:main_C}. Note that the above result only shows that the mean-square residual diverges, and not that the iterates diverge. In Example~\ref{ex:reflection}, the divergence is driven by rare trajectories along which $u_k$ grows, and the iterates may still converge almost surely.

\section{Applications}\label{sec:applications}
In this section, we discuss three applications of our results: stochastic gradient descent for convex functions, linear stochastic approximation, and learning in cocoercive games. For each application, our focus is on what the noise conditions in Assumption~\ref{assm:noise} mean, and on where they arise. Noise condition $\mathrm C$ is the most common noise model, and arises whenever the noisy evaluations are obtained from samples with bounded variance, such as stochastic gradients with bounded variance or noisy observations of a fixed quantity. In contrast, whether noise conditions $\mathrm R$ and $\mathrm I$ hold depends on the structure of the problem, and we discuss these conditions wherever applicable.

\paragraph{Root-finding formulation.} Many applications of nonexpansive SA arise in the root-finding formulation, where the goal is to find a zero of an operator $F:\mathbb{R}^d\to\mathbb{R}^d$ using noisy evaluations $F(x_k)+\zeta_{k+1}$, where $\EE[\zeta_{k+1}\mid\mathcal{F}_k]=0$, via the iteration
\begin{align}
x_{k+1}=x_k-\beta_k'\left(F(x_k)+\zeta_{k+1}\right).\label{eq:app-iteration}
\end{align}
Fix $\eta>0$, and write $\beta_k'=\eta\beta_k$. Then, \eqref{eq:app-iteration} can be rewritten as
\begin{align*}
x_{k+1}=x_k+\beta_k\left(f(x_k)-x_k+M_{k+1}\right),\qquad\text{where }f(x)=x-\eta F(x)\text{ and }M_{k+1}=-\eta\zeta_{k+1},
\end{align*}
which is iteration \eqref{eq:iteration}. The zeros of $F$ are exactly the fixed points of $f$. Moreover, the residual is $\|f(x_k)-x_k\|^2=\eta^2\|F(x_k)\|^2$, and hence our results give last-iterate bounds on $\EE\left[\|F(x_k)\|^2\right]$. The fixed-point formulation is nonexpansive when $F$ is cocoercive. An operator $F$ is said to be $(1/L)$-cocoercive if
\begin{align*}
\langle F(x)-F(y),x-y\rangle\ge\frac1L\|F(x)-F(y)\|^2,\qquad\forall x,y\in\mathbb{R}^d.
\end{align*}
If $F$ is $(1/L)$-cocoercive and has at least one zero, and $0<\eta\le2/L$, then $f$ satisfies Assumption~\ref{assm:nonexp}. Indeed, for all $x,y\in\mathbb{R}^d$,
\begin{align*}
\|f(x)-f(y)\|^2&=\|x-y\|^2-2\eta\langle F(x)-F(y),x-y\rangle+\eta^2\|F(x)-F(y)\|^2\\
&\le\|x-y\|^2-\eta\left(\frac2L-\eta\right)\|F(x)-F(y)\|^2\le\|x-y\|^2.
\end{align*}
All three applications below are of this form.

\subsection{Stochastic Gradient Descent for Convex Functions}\label{sec:app-sgd}
Consider the problem of minimizing $F_0(x)=\EE_\xi[\ell(x;\xi)]$, where $\ell$ is a loss function and $\xi$ is a random sample. Suppose that $F_0:\mathbb{R}^d\to\mathbb{R}$ is convex and $L$-smooth, and its set of minimizers is nonempty. In this case, $\nabla F_0$ is $(1/L)$-cocoercive~\cite[Theorem~2.1.5]{Nesterov2018}. Stochastic gradient descent (SGD) is given by
\begin{align*}
x_{k+1}=x_k-\beta_k'\nabla\ell(x_k;\xi_{k+1}),
\end{align*}
where $\xi_1,\xi_2,\ldots$ are i.i.d.\ samples. Since $\nabla F_0(x)=\EE_\xi[\nabla\ell(x;\xi)]$, this is of the form \eqref{eq:app-iteration} with $F=\nabla F_0$ and $\zeta_{k+1}=\nabla\ell(x_k;\xi_{k+1})-\nabla F_0(x_k)$, and our results give last-iterate bounds on $\EE\left[\|\nabla F_0(x_k)\|^2\right]$.
\begin{itemize}
\item \emph{Noise condition $\mathrm R$} is equivalent to the strong growth condition~\cite{VaswaniBachSchmidt2019},
\begin{align*}
\EE_\xi\left[\|\nabla\ell(x;\xi)\|^2\right]\le\rho\|\nabla F_0(x)\|^2\quad\text{for all } x\in\mathbb{R}^d,
\end{align*}
with $\sigma^2=\rho-1$. That is, the gradients of the individual losses are small whenever their average is small. This holds, for example, under the stronger condition $\sup_\xi\|\nabla\ell(x;\xi)\|\le B\|\nabla F_0(x)\|$ for some $B\ge0$, as previously considered in~\cite{SchmidtLeRoux2013}.
\item \emph{Noise condition $\mathrm I$} is much more general than noise condition $\mathrm R$, as it does not require the stochastic gradients to vanish at the minimizers. For example, it holds when, for some $L',G\ge0$, $\nabla\ell(\cdot\,;\xi)$ is $L'$-Lipschitz and $\|\nabla\ell(0;\xi)\|\le G$ for every $\xi$. In this case, $\|\nabla\ell(x;\xi)\|\le G+L'\|x\|$, and hence the stochastic gradients grow at most linearly in $\|x\|$.
\end{itemize}

\subsection{Linear Stochastic Approximation}\label{sec:app-linear}
Consider the problem of solving $Ax=b$, where the system has at least one solution, $A\in\mathbb{R}^{d\times d}$ and $b\in\mathbb{R}^d$, and $A$ is $(1/L)$-cocoercive, i.e., $\langle x,Ax\rangle\ge\frac1L\|Ax\|^2$ for all $x$. For example, this holds for any symmetric positive semidefinite matrix $A$ with largest eigenvalue at most $L$~\cite[p.~79]{FacchineiPang2003}. In particular, $A$ may be singular, in which case there are multiple solutions. At each time $k$, we observe conditionally unbiased estimates $\tilde A_{k+1}$ and $\tilde b_{k+1}$ of $A$ and $b$, and update
\begin{align*}
x_{k+1}=x_k+\beta_k'\left(\tilde b_{k+1}-\tilde A_{k+1}x_k\right).
\end{align*}
This is of the form \eqref{eq:app-iteration} with $F(x)=Ax-b$ and $\zeta_{k+1}=(\tilde A_{k+1}-A)x_k-(\tilde b_{k+1}-b)$, and our results give last-iterate bounds on $\EE\left[\|Ax_k-b\|^2\right]$.
Noise condition $\mathrm I$ is the standard noise condition for linear SA~\cite{SrikantYing2019,ChenMaguluriZubeldia2025}. It holds whenever $\tilde A_{k+1}$ and $\tilde b_{k+1}$ have bounded conditional second moments, since the noise term $(\tilde A_{k+1}-A)x_k$ grows linearly in $x_k$.

\subsection{Learning in Cocoercive Games}\label{sec:app-games}
Consider an $N$-player game, where each player $n\in\{1,\ldots,N\}$ chooses an action $x_n\in\mathbb{R}^d$ and incurs a cost $J_n(\mathbf{x})$. Here, $\mathbf{x}=(x_1,\ldots,x_N)\in\mathbb{R}^{Nd}$ denotes the concatenation of the actions of all players, and $J_n$ is convex and differentiable in the player's own action $x_n$. Let
\begin{align*}
F(\mathbf{x})=\left(\nabla_{x_1}J_1(\mathbf{x}),\ldots,\nabla_{x_N}J_N(\mathbf{x})\right)\in\mathbb{R}^{Nd}
\end{align*}
denote the game gradient operator. The Nash equilibria of the game are the zeros of $F$. When $F$ is $(1/L)$-cocoercive, the game is called a cocoercive game~\cite{LinEtAl2020,ChandakTamizholiBambos2026}. We assume that the game is cocoercive and has at least one Nash equilibrium. In stochastic gradient play, each player updates its action using a noisy estimate $\tilde F_{n,k+1}$ of its own gradient $\nabla_{x_n}J_n(\mathbf{x}_k)$,
\begin{align*}
x_{n,k+1}=x_{n,k}-\beta_k'\tilde F_{n,k+1},\qquad n=1,\ldots,N.
\end{align*}
Stacking the updates of all players, we get
\begin{align*}
\mathbf{x}_{k+1}=\mathbf{x}_k-\beta_k'\tilde F_{k+1},
\end{align*}
where $\mathbf{x}_k=(x_{1,k},\ldots,x_{N,k})\in\mathbb{R}^{Nd}$ and $\tilde F_{k+1}=(\tilde F_{1,k+1},\ldots,\tilde F_{N,k+1})$. This is of the form \eqref{eq:app-iteration} with $\zeta_{k+1}=\tilde F_{k+1}-F(\mathbf{x}_k)$, and our results give last-iterate bounds on $\EE\left[\|F(\mathbf{x}_k)\|^2\right]$. As the game may have multiple Nash equilibria, the residual is the natural measure of convergence in this setting.
\begin{itemize}
\item \emph{Noise condition $\mathrm R$} holds when players update at random times. Suppose each player updates at time $k$ with probability $q\in(0,1)$, independently of the other players and of the past, and scales its exact gradient by $1/q$ when it updates. Then, $\tilde F_{k+1}$ is conditionally unbiased, and $\EE\left[\|\zeta_{k+1}\|^2\mid\mathcal{F}_k\right]=(1/q-1)\|F(\mathbf{x}_k)\|^2$.
\item \emph{Noise condition $\mathrm I$} holds when the costs depend on random data, as in linear SA. For example, in quadratic games where the cost matrices are observed only through unbiased samples, the noise grows linearly in the actions.
\end{itemize}

\section{Conclusion}\label{sec:conclusion}
In this paper, we studied nonexpansive SA under the Euclidean norm. We obtained last-iterate mean-square residual bounds under three conditions on the conditional second moment of the noise. We also obtained matching lower bounds, up to logarithmic factors, and further showed that the mean-square residual may diverge under iterate-dependent noise if the stepsizes are not square-summable. Our analysis relies on a noiseless auxiliary iteration, which we use to convert an ergodic residual bound into a last-iterate bound. Future directions include obtaining optimal last-iterate rates under general norms and for two-time-scale SA with nonexpansive mappings, where it remains open whether the existing bounds~\cite{BravoCominetti2024,BlaserZhang2026,Chandak2026Nonexpansive,SarkarAggarwal2026} can be improved.

\subsection*{Disclosure of AI Use}
The authors formulated the problem and the three noise conditions, wrote preliminary versions of the proofs of the upper bounds, and conjectured that, under noise condition $\mathrm I$, square-summable stepsizes are necessary for convergence. The large language models GPT-6 Astra (OpenAI) and Claude Opus 5.5 (Anthropic) were then used to strengthen the upper bounds, and to construct the instances in Section~\ref{sec:lower}. The tools were also used to assist with writing and editing the paper. The authors have checked all results, proofs, and references, and take full responsibility for the content of this paper.

\appendix
\section{Proofs of the Upper Bounds}\label{app:restart-proofs}
Throughout this appendix, $x^\star\in\Xstar$ is fixed, and we use the notation $r_i$ and $\bar r_{\ell,k}$ introduced in Section~\ref{sec:outline}. We also repeatedly use the identity
\begin{align}
\|(1-\beta)a+\beta b\|^2=(1-\beta)\|a\|^2+\beta\|b\|^2-\beta(1-\beta)\|a-b\|^2,\qquad a,b\in\mathbb{R}^d,\ \beta\in\mathbb{R}.\label{eq:convex-identity}
\end{align}

\subsection{Proof of Lemma~\ref{lem:lyapunov}}
Subtracting $x^\star$ from both sides of \eqref{eq:iteration} and using $f(x^\star)=x^\star$, we get
\begin{align*}
x_{k+1}-x^\star=(1-\beta_k)(x_k-x^\star)+\beta_k(f(x_k)-f(x^\star))+\beta_kM_{k+1}.
\end{align*}
Applying \eqref{eq:convex-identity} with $a=x_k-x^\star$ and $b=f(x_k)-f(x^\star)$, we have
\begin{align*}
\|(1-\beta_k)(x_k-x^\star)+\beta_k(f(x_k)-f(x^\star))\|^2&=(1-\beta_k)\|x_k-x^\star\|^2+\beta_k\|f(x_k)-f(x^\star)\|^2\\
&\quad-\beta_k(1-\beta_k)\|f(x_k)-x_k\|^2\\
&\le\|x_k-x^\star\|^2-\beta_k(1-\beta_k)\|f(x_k)-x_k\|^2.
\end{align*}
Here, the inequality follows from the nonexpansiveness of $f$. Expanding the square after adding $\beta_kM_{k+1}$, we get
\begin{align*}
\|x_{k+1}-x^\star\|^2&\le\|x_k-x^\star\|^2-\beta_k(1-\beta_k)\|f(x_k)-x_k\|^2+\beta_k^2\|M_{k+1}\|^2\\
&\quad+2\beta_k\langle(1-\beta_k)(x_k-x^\star)+\beta_k(f(x_k)-x^\star),M_{k+1}\rangle.
\end{align*}
Since $x_k$ is $\mathcal{F}_k$-measurable and $\{M_{k+1}\}$ is a martingale difference sequence with respect to $\mathcal{F}_k$, the conditional expectation of the last term given $\mathcal{F}_k$ is zero. Taking expectations completes the proof of Lemma~\ref{lem:lyapunov}.

\subsection{Proof of Lemma~\ref{lem:residual-sum}}
For part (a), taking expectations in noise condition $\mathrm R$ gives $\EE\left[\|M_{i+1}\|^2\right]\le\sigma^2r_i$. Substituting this bound in Lemma~\ref{lem:lyapunov}, we get
\begin{align*}
\EE\left[\|x_{i+1}-x^\star\|^2\right]\le\EE\left[\|x_i-x^\star\|^2\right]-\beta_i[1-(1+\sigma^2)\beta_i]r_i\le\EE\left[\|x_i-x^\star\|^2\right]-[1-(1+\sigma^2)\beta_0]\beta_ir_i.
\end{align*}
Here, the second inequality follows from the fact that the stepsizes are nonincreasing. Summing from $i=0$ to $k$, the squared distance terms telescope, and we get
\begin{align*}
[1-(1+\sigma^2)\beta_0]\left(\sum_{i=0}^k\beta_i\right)\bar r_{0,k}=[1-(1+\sigma^2)\beta_0]\sum_{i=0}^k\beta_ir_i\le\|x_0-x^\star\|^2-\EE\left[\|x_{k+1}-x^\star\|^2\right]\le\|x_0-x^\star\|^2.
\end{align*}
Since $(1+\sigma^2)\beta_0<1$, dividing both sides by $[1-(1+\sigma^2)\beta_0]\sum_{i=0}^k\beta_i$ completes the proof for part (a).

For part (b), noise condition $\mathrm C$ gives $\EE\left[\|M_{i+1}\|^2\right]\le\sigma^2$. Using this bound and $1-\beta_i\ge1/2$ in Lemma~\ref{lem:lyapunov}, we get
\begin{align*}
\frac12\beta_ir_i\le\EE\left[\|x_i-x^\star\|^2\right]-\EE\left[\|x_{i+1}-x^\star\|^2\right]+\sigma^2\beta_i^2.
\end{align*}
Summing from $i=0$ to $k$ as before gives $\frac12\left(\sum_{i=0}^k\beta_i\right)\bar r_{0,k}\le\|x_0-x^\star\|^2+\sigma^2\sum_{i=0}^k\beta_i^2$, which completes the proof for part (b).

For part (c), using $\|x_i\|^2\le2\|x_i-x^\star\|^2+2\|x^\star\|^2$ in noise condition $\mathrm I$, we get
\begin{align}
\EE\left[\|M_{i+1}\|^2\right]\le\sigma^2\left(1+\EE\left[\|x_i\|^2\right]\right)\le\sigma^2\left(1+2\|x^\star\|^2+2\EE\left[\|x_i-x^\star\|^2\right]\right).\label{eq:I-noise}
\end{align}
Substituting this bound in Lemma~\ref{lem:lyapunov} and dropping the residual term, we get
\begin{align*}
\EE\left[\|x_{i+1}-x^\star\|^2\right]\le(1+2\sigma^2\beta_i^2)\EE\left[\|x_i-x^\star\|^2\right]+\sigma^2(1+2\|x^\star\|^2)\beta_i^2.
\end{align*}
Recall that $\Gamma_1=\sum_{i=0}^\infty\beta_i^2$, which is finite since $b>1/2$. Unrolling the above recursion, for all $k\ge0$,
\begin{align*}
\EE\left[\|x_k-x^\star\|^2\right]&\le\|x_0-x^\star\|^2\prod_{i=0}^{k-1}(1+2\sigma^2\beta_i^2)+\sigma^2(1+2\|x^\star\|^2)\sum_{i=0}^{k-1}\beta_i^2\prod_{j=i+1}^{k-1}(1+2\sigma^2\beta_j^2)\\
&\le e^{2\sigma^2\Gamma_1}\left(\|x_0-x^\star\|^2+\sigma^2(1+2\|x^\star\|^2)\Gamma_1\right)=:\Gamma_2.
\end{align*}
Here, the second inequality follows from $1+u\le e^u$. This proves the bound on the iterates in part (c). Substituting it in \eqref{eq:I-noise} gives $\EE\left[\|M_{i+1}\|^2\right]\le\sigma^2(1+2\|x^\star\|^2+2\Gamma_2)$ for all $i\ge0$. The rest of the proof is identical to that of part (b), with $\sigma^2$ replaced by $\sigma^2(1+2\|x^\star\|^2+2\Gamma_2)$, and using $\sum_{i=0}^k\beta_i^2\le\Gamma_1$. This completes the proof of Lemma~\ref{lem:residual-sum}.

\subsection{Proof of Lemma~\ref{lem:copy-diss}}
Throughout this proof and the next, we write $\Delta_i:=f(x_i)-x_i-f(y_i)+y_i$. Subtracting \eqref{eq:copy} from \eqref{eq:iteration}, for $k\ge\ell$,
\begin{align*}
x_{k+1}-y_{k+1}=(1-\beta_k)(x_k-y_k)+\beta_k(f(x_k)-f(y_k))+\beta_kM_{k+1}.
\end{align*}
Applying \eqref{eq:convex-identity} with $a=x_k-y_k$ and $b=f(x_k)-f(y_k)$, and using the nonexpansiveness of $f$ as in the proof of Lemma~\ref{lem:lyapunov}, we get
\begin{align*}
\|(1-\beta_k)(x_k-y_k)+\beta_k(f(x_k)-f(y_k))\|^2\le\|x_k-y_k\|^2-\beta_k(1-\beta_k)\|\Delta_k\|^2.
\end{align*}
Since $y_k$ is a deterministic function of $x_\ell$, it is $\mathcal{F}_k$-measurable. Hence, expanding $\|x_{k+1}-y_{k+1}\|^2$ and taking expectations as in the proof of Lemma~\ref{lem:lyapunov} proves \eqref{eq:copy-diss}.

Next, summing \eqref{eq:copy-diss} from $i=\ell$ to $k$, the squared distance terms telescope. Since $y_\ell=x_\ell$, we get
\begin{align*}
\sum_{i=\ell}^k\beta_i(1-\beta_i)\EE\left[\|\Delta_i\|^2\right]\le\sum_{i=\ell}^k\beta_i^2\EE\left[\|M_{i+1}\|^2\right]-\EE\left[\|x_{k+1}-y_{k+1}\|^2\right]\le\sum_{i=\ell}^k\beta_i^2\EE\left[\|M_{i+1}\|^2\right].
\end{align*}
Using $1-\beta_i\ge1/2$ proves \eqref{eq:copy-sum}.

Finally, \eqref{eq:copy} gives $y_{k+1}-y_k=\beta_k(f(y_k)-y_k)$, and hence $f(y_{k+1})-y_{k+1}=f(y_{k+1})-f(y_k)+(1-\beta_k)(f(y_k)-y_k)$. Then,
\begin{align*}
\|f(y_{k+1})-y_{k+1}\|\le\|y_{k+1}-y_k\|+(1-\beta_k)\|f(y_k)-y_k\|=\|f(y_k)-y_k\|.
\end{align*}
Here, the inequality follows from the triangle inequality and the nonexpansiveness of $f$. Iterating from $k=\ell$ and using $y_\ell=x_\ell$ completes the proof of Lemma~\ref{lem:copy-diss}.

\subsection{Proof of Lemma~\ref{lem:copy-future}}
Fix $0\le\ell<k$, and let $\{y_i\}_{i\ge\ell}$ be the auxiliary iteration \eqref{eq:copy} restarted at time $\ell$. For $i\ge\ell$, note that $f(y_i)-y_i=(f(x_i)-x_i)-\Delta_i$. Then, by Lemma~\ref{lem:copy-diss},
\begin{align*}
\|f(x_\ell)-x_\ell\|^2\ge\|f(y_i)-y_i\|^2&=\|f(x_i)-x_i\|^2-2\langle f(x_i)-x_i,\Delta_i\rangle+\|\Delta_i\|^2\\
&\ge\|f(x_i)-x_i\|^2-2\langle f(x_i)-x_i,\Delta_i\rangle.
\end{align*}
Taking expectations, multiplying by $\beta_i$, and summing from $i=\ell$ to $k$, we get
\begin{align*}
\left(\sum_{i=\ell}^k\beta_i\right)\left(\bar r_{\ell,k}-r_\ell\right)&\le2\sum_{i=\ell}^k\beta_i\EE\left[\langle f(x_i)-x_i,\Delta_i\rangle\right]\le2\left(\sum_{i=\ell}^k\beta_ir_i\right)^{1/2}\left(\sum_{i=\ell}^k\beta_i\EE\left[\|\Delta_i\|^2\right]\right)^{1/2}\\
&\le2\sqrt2\left(\sum_{i=\ell}^k\beta_i\right)^{1/2}\left(\bar r_{\ell,k}\right)^{1/2}\left(\sum_{i=\ell}^k\beta_i^2\EE\left[\|M_{i+1}\|^2\right]\right)^{1/2}.
\end{align*}
Here, the second inequality follows from the Cauchy--Schwarz inequality, first for the expectation and then for the weighted sum, and the third inequality follows from \eqref{eq:copy-sum}. Dividing by $\sum_{i=\ell}^k\beta_i$, we get
\begin{align}
\bar r_{\ell,k}-r_\ell\le2\sqrt2\left(\bar r_{\ell,k}\right)^{1/2}\left(\frac{\sum_{i=\ell}^k\beta_i^2\EE\left[\|M_{i+1}\|^2\right]}{\sum_{i=\ell}^k\beta_i}\right)^{1/2}.\label{eq:future-excess}
\end{align}
Next, splitting off the first term of the weighted sum gives $\left(\sum_{i=\ell}^k\beta_i\right)\bar r_{\ell,k}=\beta_\ell r_\ell+\left(\sum_{i=\ell+1}^k\beta_i\right)\bar r_{\ell+1,k}$. Writing $\sum_{i=\ell}^k\beta_i=\beta_\ell+\sum_{i=\ell+1}^k\beta_i$ and rearranging, we get
\begin{align*}
\bar r_{\ell+1,k}=\bar r_{\ell,k}+\frac{\beta_\ell}{\sum_{i=\ell+1}^k\beta_i}\left(\bar r_{\ell,k}-r_\ell\right).
\end{align*}
Substituting \eqref{eq:future-excess} completes the proof of Lemma~\ref{lem:copy-future}.

\subsection{Proof of Lemma~\ref{lem:copy-propagate}}
Fix $0\le m\le k$. When $m=k$, we have $r_k=\bar r_{k,k}$, and hence part (a) holds with equality, and part (b) holds for any $\Gamma_3\ge1$. We henceforth assume that $m<k$, and write
\begin{align*}
c_i:=\frac{\beta_i\sqrt{\beta_i}}{\sum_{j=i+1}^k\beta_j},\quad m\le i<k,\qquad\text{and}\qquad S:=\sum_{i=m}^{k-1}c_i>0.
\end{align*}
For part (a), taking expectations in noise condition $\mathrm R$ gives $\EE\left[\|M_{j+1}\|^2\right]\le\sigma^2r_j$. Since the stepsizes are nonincreasing, for $m\le i<k$,
\begin{align*}
\sum_{j=i}^k\beta_j^2\EE\left[\|M_{j+1}\|^2\right]\le\sigma^2\beta_i\sum_{j=i}^k\beta_jr_j=\sigma^2\beta_i\left(\sum_{j=i}^k\beta_j\right)\bar r_{i,k}.
\end{align*}
Substituting this bound in Lemma~\ref{lem:copy-future} with $\ell=i$, we get $\bar r_{i+1,k}\le\left(1+2\sqrt2\sigma c_i\right)\bar r_{i,k}$ for $m\le i<k$. Iterating this inequality from $i=m$ to $k-1$ and using $\bar r_{k,k}=r_k$ completes the proof for part (a).

For part (b), suppose that $\EE\left[\|M_{i+1}\|^2\right]\le\bar\sigma^2$ for all $i\ge0$, for some $\bar\sigma>0$. Under noise condition $\mathrm C$, this holds with $\bar\sigma^2=\sigma^2$. Under noise condition $\mathrm I$ with $1/2<b<1$, Lemma~\ref{lem:residual-sum}(c) and \eqref{eq:I-noise} show that this holds with $\bar\sigma^2=\sigma^2(1+2\|x^\star\|^2+2\Gamma_2)$. Since the stepsizes are nonincreasing, $\sum_{j=i}^k\beta_j^2\EE\left[\|M_{j+1}\|^2\right]\le\bar\sigma^2\beta_i\sum_{j=i}^k\beta_j$. Substituting this bound in Lemma~\ref{lem:copy-future} with $\ell=i$, and using the AM--GM inequality, we get
\begin{align*}
\bar r_{i+1,k}\le\bar r_{i,k}+2\sqrt2\bar\sigma c_i\left(\bar r_{i,k}\right)^{1/2}\le\left(1+\frac{c_i}{S}\right)\bar r_{i,k}+2\bar\sigma^2c_iS,\qquad m\le i<k.
\end{align*}
Iterating this inequality from $i=m$ to $k-1$, and bounding the products multiplying the additive terms by the full product, as every multiplicative factor is at least one, we get
\begin{align*}
r_k=\bar r_{k,k}\le\prod_{i=m}^{k-1}\left(1+\frac{c_i}{S}\right)\left(\bar r_{m,k}+2\bar\sigma^2S\sum_{i=m}^{k-1}c_i\right)\le e\left(\bar r_{m,k}+2\bar\sigma^2S^2\right).
\end{align*}
Here, the second inequality follows from $1+u\le e^u$ and $\sum_{i=m}^{k-1}c_i=S$. This completes the proof for part (b) with $\Gamma_3=e\max\{1,2\bar\sigma^2\}$.

\subsection{Proof of Lemma~\ref{lem:terminal-sum}}
Recall that $m=\lceil k/2\rceil$ and $\beta_i=\beta/(i+K_0)^b$. Since $r_i\ge0$ for all $i$,
\begin{align*}
\bar r_{m,k}=\frac{\sum_{i=m}^k\beta_ir_i}{\sum_{i=m}^k\beta_i}\le\frac{\sum_{i=0}^k\beta_ir_i}{\sum_{i=m}^k\beta_i}=\frac{\sum_{i=0}^k\beta_i}{\sum_{i=m}^k\beta_i}\,\bar r_{0,k}.
\end{align*}
Hence, it suffices to bound the ratio of the stepsize sums. Since $k-m+1\ge(k+1)/2$ and the stepsizes are nonincreasing,
\begin{align*}
\sum_{i=m}^k\beta_i\ge(k-m+1)\beta_k\ge\frac{(k+1)\beta}{2(k+K_0)^b}.
\end{align*}
Moreover, comparing the sum with an integral, we get
\begin{align*}
\sum_{i=0}^k\beta_i\le\frac{\beta}{K_0^b}+\beta\int_{K_0}^{k+K_0}s^{-b}\,ds\le\frac{\beta}{K_0^b}+\frac{\beta(k+K_0)^{1-b}}{1-b}.
\end{align*}
Combining the above two bounds and using $k+K_0\le K_0(k+1)$, which holds since $K_0\ge1$, we get
\begin{align*}
\frac{\sum_{i=0}^k\beta_i}{\sum_{i=m}^k\beta_i}\le\frac{2}{k+1}\left(\frac{k+K_0}{K_0}\right)^b+\frac{2(k+K_0)}{(1-b)(k+1)}\le2+\frac{2K_0}{1-b}=:\Gamma_4.
\end{align*}
This proves the first inequality in \eqref{eq:terminal}.

For the second inequality, the sum is empty for $k\le1$, and we assume that $k\ge2$. For $m\le i<k$, note that $\beta_i/\beta_k=((k+K_0)/(i+K_0))^b\le2^b$ and $\sum_{j=i+1}^k\beta_j\ge(k-i)\beta_k$. Therefore,
\begin{align*}
\sum_{i=m}^{k-1}\frac{\beta_i\sqrt{\beta_i}}{\sum_{j=i+1}^k\beta_j}\le2^{3b/2}\sqrt{\beta_k}\sum_{i=m}^{k-1}\frac{1}{k-i}\le2^{3b/2}\sqrt{\beta_k}(1+\log(k-m))\le\frac{3\cdot2^{3b/2}\sqrt\beta\log(k+2)}{(k+1)^{b/2}}.
\end{align*}
Here, the second inequality follows by comparing the harmonic sum with an integral, and the third inequality follows from $\sqrt{\beta_k}\le\sqrt\beta/(k+1)^{b/2}$ and $1+\log(k-m)\le3\log(k+2)$. This completes the proof of Lemma~\ref{lem:terminal-sum} with $\Gamma_5=3\cdot2^{3b/2}\sqrt\beta$.

\subsection{Proof of Theorem~\ref{thm:main_R}}
Throughout the proofs of Theorems~\ref{thm:main_R}--\ref{thm:main_I}, we take $m=\lceil k/2\rceil$. Since $K_0\ge1$, we have $(k+K_0)^b\le K_0^b(k+1)^b$, and hence $\sum_{i=0}^k\beta_i\ge(k+1)\beta_k\ge\beta(k+1)^{1-b}/K_0^b$. Combining this with the first inequality in Lemma~\ref{lem:terminal-sum}, we get
\begin{align}
\bar r_{m,k}\le\Gamma_4\bar r_{0,k}\le\frac{\Gamma_6}{(k+1)^{1-b}}\left(\sum_{i=0}^k\beta_i\right)\bar r_{0,k},\qquad\text{where }\Gamma_6:=\frac{\Gamma_4K_0^b}{\beta}.\label{eq:rbar-m}
\end{align}
Under noise condition $\mathrm R$, Lemma~\ref{lem:residual-sum}(a) then gives $\bar r_{m,k}\le\Gamma_6\|x_0-x^\star\|^2/([1-(1+\sigma^2)\beta_0](k+1)^{1-b})$. For the product in Lemma~\ref{lem:copy-propagate}(a), using $1+u\le e^u$ and the second inequality in Lemma~\ref{lem:terminal-sum}, we get
\begin{align*}
\prod_{i=m}^{k-1}\left[1+2\sqrt2\sigma\frac{\beta_i\sqrt{\beta_i}}{\sum_{j=i+1}^k\beta_j}\right]\le\exp\left(\frac{2\sqrt2\sigma\Gamma_5\log(k+2)}{(k+1)^{b/2}}\right)\le\exp\left(2\sqrt2\sigma\Gamma_5\left(\log2+\frac{2}{be}\right)\right).
\end{align*}
Here, the second inequality follows from $\log(k+2)\le\log2+\log(k+1)$ and $\sup_{t\ge1}\log t/t^{b/2}=2/(be)$. Substituting the above two bounds in Lemma~\ref{lem:copy-propagate}(a), we get
\begin{align*}
r_k\le\frac{\Gamma_6\|x_0-x^\star\|^2}{1-(1+\sigma^2)\beta_0}\exp\left(2\sqrt2\sigma\Gamma_5\left(\log2+\frac{2}{be}\right)\right)\frac{1}{(k+1)^{1-b}}.
\end{align*}
This completes the proof of the upper bound in Theorem~\ref{thm:main_R} with
\begin{align*}
\cfrak_1=\frac{\Gamma_6\|x_0-x^\star\|^2}{1-(1+\sigma^2)\beta_0}\exp\left(2\sqrt2\sigma\Gamma_5\left(\log2+\frac{2}{be}\right)\right).
\end{align*}
The optimality of this bound follows from Proposition~\ref{prop:copy-lower-R}.

\subsection{Proof of Theorem~\ref{thm:main_C}}
We first bound the sum of squared stepsizes. Since $K_0\ge1$, comparing the sum with an integral gives
\begin{align*}
\sum_{i=0}^k\beta_i^2\le\beta^2\left(1+\int_1^{k+1}s^{-2b}\,ds\right)\le\begin{cases}
\beta^2(k+1)^{1-2b}/(1-2b),&0<b<\tfrac12,\\
3\beta^2\log(k+2),&b=\tfrac12,\\
\Gamma_1,&\tfrac12<b<1,
\end{cases}
\end{align*}
where $\Gamma_1=\sum_{i=0}^\infty\beta_i^2$, and the second case uses $1+\log(k+1)\le3\log(k+2)$. Using Lemma~\ref{lem:copy-propagate}(b) with $\Gamma_3=e\max\{1,2\sigma^2\}$, \eqref{eq:rbar-m}, Lemma~\ref{lem:residual-sum}(b), and the second inequality in Lemma~\ref{lem:terminal-sum}, we get
\begin{align}
r_k\le\Gamma_3\left(\frac{\Gamma_6}{(k+1)^{1-b}}\left(2\|x_0-x^\star\|^2+2\sigma^2\sum_{i=0}^k\beta_i^2\right)+\frac{\Gamma_5^2\log^2(k+2)}{(k+1)^b}\right).\label{eq:C-combined}
\end{align}
For $0<b<1/2$, using the first case of the bound on $\sum_{i=0}^k\beta_i^2$ and $(k+1)^{1-b}\ge(k+1)^b$ in \eqref{eq:C-combined}, we get
\begin{align*}
r_k&\le\frac{\Gamma_3}{(k+1)^b}\left(\Gamma_6\left(2\|x_0-x^\star\|^2+\frac{2\sigma^2\beta^2}{1-2b}\right)+\Gamma_5^2\log^2(k+2)\right)\\
&\le\Gamma_3\left(\frac{\Gamma_6}{\log^22}\left(2\|x_0-x^\star\|^2+\frac{2\sigma^2\beta^2}{1-2b}\right)+\Gamma_5^2\right)\frac{\log^2(k+2)}{(k+1)^b}.
\end{align*}
For $b=1/2$, using the second case of the bound on $\sum_{i=0}^k\beta_i^2$ in \eqref{eq:C-combined}, we similarly get
\begin{align*}
r_k\le\Gamma_3\left(\frac{2\Gamma_6\|x_0-x^\star\|^2}{\log^22}+\frac{6\Gamma_6\sigma^2\beta^2}{\log2}+\Gamma_5^2\right)\frac{\log^2(k+2)}{(k+1)^{1/2}}.
\end{align*}
This proves the first case of Theorem~\ref{thm:main_C} with
\begin{align*}
\cfrak_2=\begin{cases}
\Gamma_3\left(\dfrac{\Gamma_6}{\log^22}\left(2\|x_0-x^\star\|^2+\dfrac{2\sigma^2\beta^2}{1-2b}\right)+\Gamma_5^2\right),&0<b<\tfrac12,\\[8pt]
\Gamma_3\left(\dfrac{2\Gamma_6\|x_0-x^\star\|^2}{\log^22}+\dfrac{6\Gamma_6\sigma^2\beta^2}{\log2}+\Gamma_5^2\right),&b=\tfrac12.
\end{cases}
\end{align*}
For $1/2<b<1$, we first note that
\begin{align}
\frac{\log^2(k+2)}{(k+1)^{2b-1}}\le\frac{2\log^22+2\log^2(k+1)}{(k+1)^{2b-1}}\le2\log^22+\frac{8}{e^2(2b-1)^2}.\label{eq:log-ratio}
\end{align}
Here, the first inequality follows from $\log(k+2)\le\log2+\log(k+1)$, and the second inequality follows from $\sup_{t\ge1}\log^2t/t^{2b-1}=4/(e^2(2b-1)^2)$. Using the third case of the bound on $\sum_{i=0}^k\beta_i^2$ and \eqref{eq:log-ratio} in \eqref{eq:C-combined}, we get
\begin{align*}
r_k\le\Gamma_3\left(\Gamma_6(2\|x_0-x^\star\|^2+2\sigma^2\Gamma_1)+\Gamma_5^2\left(2\log^22+\frac{8}{e^2(2b-1)^2}\right)\right)\frac{1}{(k+1)^{1-b}}.
\end{align*}
This completes the proof of the upper bound in Theorem~\ref{thm:main_C} with
\begin{align*}
\cfrak_3=\Gamma_3\left(\Gamma_6(2\|x_0-x^\star\|^2+2\sigma^2\Gamma_1)+\Gamma_5^2\left(2\log^22+\frac{8}{e^2(2b-1)^2}\right)\right).
\end{align*}
The optimality of this bound, up to a logarithmic factor, follows from Propositions~\ref{prop:copy-lower-R} and~\ref{prop:copy-lower-C}.

\subsection{Proof of Theorem~\ref{thm:main_I}}
The proof is identical to that of the case $1/2<b<1$ in Theorem~\ref{thm:main_C}, with Lemma~\ref{lem:residual-sum}(c) in place of Lemma~\ref{lem:residual-sum}(b), and with $\Gamma_3=e\max\{1,2\sigma^2(1+2\|x^\star\|^2+2\Gamma_2)\}$ in Lemma~\ref{lem:copy-propagate}(b). Using \eqref{eq:rbar-m}, Lemma~\ref{lem:residual-sum}(c), the second inequality in Lemma~\ref{lem:terminal-sum}, and \eqref{eq:log-ratio}, we get
\begin{align*}
r_k\le\Gamma_3\left(2\Gamma_6\left(\|x_0-x^\star\|^2+\sigma^2\Gamma_1(1+2\|x^\star\|^2+2\Gamma_2)\right)+\Gamma_5^2\left(2\log^22+\frac{8}{e^2(2b-1)^2}\right)\right)\frac{1}{(k+1)^{1-b}}.
\end{align*}
This completes the proof of the upper bound in Theorem~\ref{thm:main_I} with
\begin{align*}
\cfrak_4=\Gamma_3\left(2\Gamma_6\left(\|x_0-x^\star\|^2+\sigma^2\Gamma_1(1+2\|x^\star\|^2+2\Gamma_2)\right)+\Gamma_5^2\left(2\log^22+\frac{8}{e^2(2b-1)^2}\right)\right).
\end{align*}
The optimality of this bound follows from Proposition~\ref{prop:copy-lower-R}, and the divergence for $0<b<1/2$ follows from Proposition~\ref{prop:copy-lower-I}.

\section{Proofs of the Optimality Results}\label{app:lower-proofs}
Recall that $r_i=\EE\left[\|f(x_i)-x_i\|^2\right]$ denotes the mean-square residual at time $i$. Throughout this appendix, we use the facts about the rotation $R_\theta$ stated at the beginning of Section~\ref{sec:lower}.

\subsection{Proof of Proposition~\ref{prop:copy-lower-R}}
Fix $k\ge1$, consider the instance in Example~\ref{ex:rotation}, and write $\gamma:=2(1-\cos\theta)=1/(1+2\sum_{j=0}^{k-1}\beta_j)$. Since $M_{i+1}=0$, the iteration is $x_{i+1}=(1-\beta_i)x_i+\beta_iR_\theta x_i$. Applying \eqref{eq:copy-single-identity} repeatedly, we get
\begin{align*}
r_k=\gamma\|x_k\|^2=\gamma\prod_{j=0}^{k-1}\left[1-\gamma\beta_j(1-\beta_j)\right].
\end{align*}
Since $\gamma\le1$ and $\beta_j\le1/2$, we have $0\le\gamma\beta_j(1-\beta_j)\le1/2$. Using $1-u\ge e^{-2u}$ for $0\le u\le1/2$, we get
\begin{align}
\prod_{j=0}^{k-1}\left[1-\gamma\beta_j(1-\beta_j)\right]\ge\exp\left(-2\gamma\sum_{j=0}^{k-1}\beta_j\right)\ge e^{-1}.\label{eq:copy-single-product}
\end{align}
Here, the second inequality follows from $2\gamma\sum_{j=0}^{k-1}\beta_j=2\sum_{j=0}^{k-1}\beta_j/(1+2\sum_{j=0}^{k-1}\beta_j)\le1$. Moreover, since $K_0\ge1$, comparing the sum with an integral gives
\begin{align*}
\sum_{j=0}^{k-1}\beta_j\le\beta\sum_{j=0}^{k-1}(j+1)^{-b}\le\beta\left(1+\int_1^ks^{-b}\,ds\right)\le\frac{\beta(k+1)^{1-b}}{1-b},
\end{align*}
and hence $\gamma\ge\Gamma_7/(k+1)^{1-b}$, where $\Gamma_7:=(1+2\beta/(1-b))^{-1}$. Combining the above bounds gives $r_k\ge\Gamma_7/(e(k+1)^{1-b})$. This completes the proof of Proposition~\ref{prop:copy-lower-R} with $\cfrak_5=\Gamma_7/e$.

\subsection{Proof of Proposition~\ref{prop:copy-lower-C}}
Fix $k\ge1$, consider the instance in Example~\ref{ex:multiscale}, and write $\gamma_p:=2(1-\cos\theta_p)$. Since $\varepsilon_{i+1}$ is independent of $\mathcal{F}_i$, we have $\EE\left[M_{i+1}\mid\mathcal{F}_i\right]=0$ and $\EE\left[\|M_{i+1}\|^2\mid\mathcal{F}_i\right]=\sigma^2$, and hence noise condition $\mathrm C$ holds. Moreover, the noise injected at times $i\ge k$ does not affect $x_k$. For $1\le p\le k$, let $x_j^{(p)}\in\mathbb{R}^2$ denote the $p$-th pair of coordinates of $x_j$.

We first compute the contribution of each pair to the residual at time $k$. Fix $1\le p\le k$. Since $x_0=0$, and before time $k$ the $p$-th pair receives noise only at time $p-1$, we have $x_j^{(p)}=0$ for $j\le p-1$, $x_p^{(p)}=\sigma\beta_{p-1}\varepsilon_p(1,0)$, and $x_{j+1}^{(p)}=(1-\beta_j)x_j^{(p)}+\beta_jR_{\theta_p}x_j^{(p)}$ for $p\le j\le k-1$. Since $\varepsilon_p^2=1$, applying \eqref{eq:copy-single-identity} repeatedly gives
\begin{align*}
\left\|R_{\theta_p}x_k^{(p)}-x_k^{(p)}\right\|^2=\gamma_p\sigma^2\beta_{p-1}^2\prod_{j=p}^{k-1}\left[1-\gamma_p\beta_j(1-\beta_j)\right]\ge\frac{\gamma_p\sigma^2\beta_{p-1}^2}{e}.
\end{align*}
Here, the inequality follows as in \eqref{eq:copy-single-product}, since $2\gamma_p\sum_{j=p}^{k-1}\beta_j\le1$. Since $f$ acts separately on each pair and the pairs are orthogonal, $\|f(x_k)-x_k\|^2=\sum_{p=1}^k\|R_{\theta_p}x_k^{(p)}-x_k^{(p)}\|^2$. Summing the above bound over the pairs and substituting $i=p-1$, we get
\begin{align}
r_k\ge\frac{\sigma^2}{e}\sum_{i=0}^{k-1}\frac{\beta_i^2}{1+2\sum_{j=i+1}^{k-1}\beta_j}.\label{eq:copy-lower-multiscale}
\end{align}
It remains to lower bound the sum in \eqref{eq:copy-lower-multiscale}. For $k=1$, we have $\gamma_1=1$, and the above computation gives $r_1=\sigma^2\beta_0^2$ exactly. For $k\ge2$, let $m=\lceil k/2\rceil$. Since all terms in \eqref{eq:copy-lower-multiscale} are nonnegative, we can restrict the sum to $m\le i<k$. For $m\le i<k$, we have $\beta_i\ge\beta_k$ and $\sum_{j=i+1}^{k-1}\beta_j\le(k-i-1)\beta_m$. Hence,
\begin{align*}
r_k\ge\frac{\sigma^2\beta_k^2}{e}\sum_{j=0}^{k-m-1}\frac{1}{1+2j\beta_m}\ge\frac{\sigma^2\beta_k^2}{e}\int_0^{k-m}\frac{ds}{1+2s\beta_m}=\frac{\sigma^2\beta_k^2}{2e\beta_m}\log\left(1+2(k-m)\beta_m\right).
\end{align*}
We now bound the two factors on the right-hand side. Since $m+K_0\ge(k+K_0)/2$ and $k+K_0\le K_0(k+1)$,
\begin{align*}
\frac{\beta_k^2}{\beta_m}=\frac{\beta(m+K_0)^b}{(k+K_0)^{2b}}\ge\frac{\beta}{2^b(k+K_0)^b}\ge\frac{\Gamma_8}{(k+1)^b},\qquad\text{where }\Gamma_8:=\frac{\beta}{2^bK_0^b}.
\end{align*}
Moreover, since $k-m\ge(k+1)/4$ for $k\ge2$, we have $2(k-m)\beta_m\ge(k+1)\beta_k/2\ge\beta(k+1)^{1-b}/(2K_0^b)$. Using $\log(1+cx)\ge\min\{1,c\}\log(1+x)$ for $c,x\ge0$, we get
\begin{align*}
\log\left(1+2(k-m)\beta_m\right)\ge\min\left\{1,\frac{\beta}{2K_0^b}\right\}\log\left(1+(k+1)^{1-b}\right)\ge\Gamma_9\log(k+1),
\end{align*}
where $\Gamma_9:=(1-b)\min\{1,\beta/(2K_0^b)\}$. Combining the above bounds gives $r_k\ge\sigma^2\Gamma_8\Gamma_9\log(k+1)/(2e(k+1)^b)$ for $k\ge2$. Together with the case $k=1$, this completes the proof of Proposition~\ref{prop:copy-lower-C} with
\begin{align*}
\cfrak_6=\min\left\{\frac{\Gamma_8\Gamma_9}{2e},\frac{2^b\beta_0^2}{\log2}\right\}.
\end{align*}

\subsection{Proof of Proposition~\ref{prop:copy-lower-I}}
Consider the instance in Example~\ref{ex:reflection}. Since $f(x)-x=(0,-2v)$, the coordinates of the iteration \eqref{eq:iteration} evolve as
\begin{align*}
u_{i+1}=u_i+\beta_iM_{i+1}^{(1)},\qquad v_{i+1}=(1-2\beta_i)v_i+\beta_iM_{i+1}^{(2)}.
\end{align*}
Since the coordinates of $M_{i+1}$ have conditional mean zero and conditional second moment $\sigma^2u_i^2/2$, taking expectations gives
\begin{align*}
\EE\left[u_{i+1}^2\right]=\left(1+\frac{\sigma^2\beta_i^2}{2}\right)\EE\left[u_i^2\right],\qquad\EE\left[v_{i+1}^2\right]=(1-2\beta_i)^2\EE\left[v_i^2\right]+\frac{\sigma^2\beta_i^2}{2}\EE\left[u_i^2\right].
\end{align*}
Since $u_0=1$, iterating the first recursion and using $\log(1+z)\ge z/(1+z)$ with $z=\sigma^2\beta_j^2/2\le\sigma^2/8$, we get
\begin{align*}
\EE\left[u_i^2\right]=\prod_{j=0}^{i-1}\left(1+\frac{\sigma^2\beta_j^2}{2}\right)\ge\exp\left(\frac{\sigma^2/2}{1+\sigma^2/8}\sum_{j=0}^{i-1}\beta_j^2\right).
\end{align*}
Dropping the nonnegative first term in the recursion for $\EE[v_k^2]$, for every $k\ge1$,
\begin{align}
r_k=4\EE\left[v_k^2\right]\ge2\sigma^2\beta_{k-1}^2\EE\left[u_{k-1}^2\right]\ge2\sigma^2\beta_{k-1}^2\exp\left(\frac{\sigma^2/2}{1+\sigma^2/8}\sum_{j=0}^{k-2}\beta_j^2\right).\label{eq:reflection-finite}
\end{align}
\textit{Part (a).} We bound the two factors on the right-hand side of \eqref{eq:reflection-finite}. Since $k-1+K_0\le K_0(k+1)$, we have $\beta_{k-1}^2\ge\Gamma_{10}/(k+1)^{2b}$, where $\Gamma_{10}:=\beta^2/K_0^{2b}$. For $k\ge2$, the monotonicity of the stepsizes and $k-1\ge(k+1)/3$ give
\begin{align*}
\sum_{j=0}^{k-2}\beta_j^2\ge(k-1)\beta_{k-2}^2\ge\frac{(k-1)\beta^2}{K_0^{2b}(k+1)^{2b}}\ge\Gamma_{11}(k+1)^{1-2b},\qquad\text{where }\Gamma_{11}:=\frac{\beta^2}{3K_0^{2b}}.
\end{align*}
Since the sum is empty for $k=1$, we have $\sum_{j=0}^{k-2}\beta_j^2\ge\Gamma_{11}((k+1)^{1-2b}-2^{1-2b})$ for all $k\ge1$. Substituting the above bounds in \eqref{eq:reflection-finite} completes the proof of Proposition~\ref{prop:copy-lower-I} with
\begin{align*}
\cfrak_8=\frac{\sigma^2\Gamma_{11}}{2(1+\sigma^2/8)},\qquad\cfrak_7=2\sigma^2\Gamma_{10}\exp\left(-2^{1-2b}\cfrak_8\right).
\end{align*}
Finally, since $1-2b>0$, the exponential term dominates the polynomial term, and hence $r_k\to\infty$ as $k\to\infty$.

\textit{Part (b).} Let $b=1/2$ and $p:=\sigma^2\beta^2/2>1$. Then $\beta_j^2=\beta^2/(j+K_0)$. Using $\log(1+z)\ge z-z^2/2$ for $z\ge0$, we get
\begin{align*}
\EE\left[u_i^2\right]=\prod_{j=0}^{i-1}\left(1+\frac{p}{j+K_0}\right)\ge\exp\left(p\sum_{j=0}^{i-1}\frac{1}{j+K_0}-\frac{p^2}{2}\sum_{j=0}^{\infty}\frac{1}{(j+K_0)^2}\right)\ge e^{-p^2}\left(\frac{i+K_0}{K_0}\right)^{p},
\end{align*}
where we used $\sum_{j=0}^{i-1}1/(j+K_0)\ge\log((i+K_0)/K_0)$ and $\sum_{j=0}^\infty1/(j+K_0)^2\le1/K_0^2+1/K_0\le2$. Substituting this with $i=k-1$ in the first inequality in \eqref{eq:reflection-finite}, and using $\beta_{k-1}^2=\beta^2/(k-1+K_0)$ and $k-1+K_0\ge(k+1)/2$, we get
\begin{align*}
r_k\ge2\sigma^2\beta^2e^{-p^2}K_0^{-p}(k-1+K_0)^{p-1}\ge\cfrak_9(k+1)^{p-1},\qquad\text{where }\cfrak_9:=2^{2-p}\sigma^2\beta^2e^{-p^2}K_0^{-p}.
\end{align*}
Since $p>1$, we get $r_k\to\infty$ as $k\to\infty$.

 \bibliographystyle{plainnat}
 \bibliography{refs}

\end{document}